\documentclass[12pt,twoside]{article}
\usepackage{smart}
\usepackage{oldamsf}
\Aiv
\plainpage
\theorems
\UseAMSsymbols
\loadscript
\userlabels=[,${}^\circ$,${}^{\circ\circ}$,${}^{\circ\circ\circ}$];%
\let\text\mbox
\usepackage{ergomath}
\usepackage{diagr2xy}
\warndef\eop{\hbox{\vrule width 6pt height 6pt depth 0pt}}
\warndef\arr{\longrightarrow}
\mathdef\sar#1,#2,#3;$#1\stackrel{#2}{\arr}#3$
\mathdef\darr#1#2#3#4$#1,#2\colon#3{}_{\arr}^{\arr}#4$
\mathdef\spar#1.#2.#3;$#1\stackrel{#2}{\arr}#3$
\mathdef\fam#1#2$\{#1\}_{#2}$
\warndef\comp{\mbox{\scriptsize$\circ$}}
\warndef\mult{\centerdot}
\mathdef\N$\Bbb N$
\mathdef\Z$\Bbb Z$
\mathdef\admon${\bf AdMon}$
\mathdef\mon${\bf Mon}$
\mathdef\admonl${\bf AdMon_L}$
\mathdef\admonb${\bf AdMonB}$
\mathdef\free${\script F(\emptyset)}$
\mathdef\freemon${\script F_{\bf Mon}(\eta_0,\eta_1,\dots,\eta_k,\dots
\varepsilon_0,\varepsilon_1,\dots,\varepsilon_k,\dots)}$
\begin{document}
\address{Institute for Nuclear Research and Nuclear Energy,\\
Tsarigradsko chosse blvd.\ 72, BG-1784 Sofia, Bulgaria.\\
e-mail: {\sl vmolot@inrne.bas.bg}}
\author{Vladimir MOLOTKOV%
\thanks{%
Financial support of the Ministry of Sciences and Education of Bulgaria
under grant F-610/96-97 and MPIMiS, Leipzig, where the final editing was performed,
is most gratefully acknowledged.}%
}
\date{}
\title{Adjunctions in Monoids}
\maketitle

Let $M$ and $N$ be monoids considered as categories with the only object.

Let
\begin{equation}
\label{1}
f\colon M\arr N,\ \ g\colon N\arr M
\end{equation}
be morphisms of monoids, considered as functors. Let the
functor $f$ is left adjoint to the functor $g$.

{\sl Is it true then that $f$ (or, what is the same, $g$) is always an
isomorphism\/}?

In \cite{M}, p.136, this question was posed as an open question.
Here I answer this question and the answer is {\sl no}.%
\footnote{This result was obtained actually somewhere at beginning of $90^{\rm th}$;
at the end of 1999 I asked Prof.~Manes in an e-mail about the status of the question.
He answered that ``to his knowledge the question is still open''.
}

To prove this, I will construct a Birkhoff variety of algebras, which
is naturally equivalent to the category of adjunctions in monoids, and
consider its initial object which is a monoid generated by
$2\times\N$ free variables subject to a certain set of relations.  An
application of M.~H.~A.~Newman's reduction theorem (\cite{N}, cited by
\cite{C}) permits one to describe the canonical form of elements in
the monoid and, in particular, to negatively answer the question
posed.

Let
\begin{equation}
\label{2}
\varphi\colon N\arr M,\ \
\end{equation}
be an isomorphism of {\it sets} such that the triple
$(f,g,\varphi)\colon M\rightharpoonup N$
is an adjunction, which means that the identity

\begin{equation}
\label{3}
g(n)\varphi(n')m = \varphi(nn'f(m))
\end{equation}
holds for every $n,n'\in N$ and $m\in M$ (see, e.g., \cite{ML}, p.78).

The identity (\ref{3}) together with the fact that $\varphi$ is an iso
implies, evidently, that

\begin{userlabel}
\addtocounter{equation}{-2}
\addtocounter{labelcounter}{1}
\begin{equation}
\label{2a}
\varphi^{-1}\colon M\arr N
\end{equation}
\end{userlabel}
satisfies the ``dual'' identity:

\begin{userlabel}
\begin{equation}
\addtocounter{labelcounter}{1}
\label{3a}
n\varphi^{-1}(m')f(m) = \varphi^{-1}(g(n)m'm)
\end{equation}
\end{userlabel}
(to get from some identity its `` dual'' replace $n$'s with $m$'s
and vice versa; $f\leftrightarrow g$, $\varphi\leftrightarrow \varphi^{-1}$
and, finally, invert the order of all compositions).

Let, further, $\eta\in M$ and $\varepsilon\in N$ be defined as

\begin{equation}
\label{4}
\eta := \varphi(1),\ \ \varepsilon := \varphi^{-1}(1).
\end{equation}
Then by setting $n'=m=1$ in (\ref{3}) we get:

\begin{userlabel}
\begin{equation}
\label{5}
\varphi(n) = g(n)\eta
\end{equation}
and, dually:
\begin{equation}
\label{5a}
\varphi^{-1}(m) = \varepsilon f(m)
\end{equation}
\end{userlabel}
i.e., $\varphi$ and $\varphi^{-1}$ are factorized as follows:
\begin{userlabel}
\begin{equation}
\label{6}
\varphi = {\rm R}_\eta\comp g\ \ \
(\varphi\colon\sar N,g,M;\sar {},{\rm R}_\eta,M;)
\end{equation}
\begin{equation}
\label{6a}
\varphi^{-1} = {\rm L}_\varepsilon\comp f\ \ \
(\varphi^{-1}\colon\sar M,f,N;\sar{},{\rm L}_\varepsilon,N;),
\end{equation}
\end{userlabel}
where ${\rm R}_\eta$ (resp. ${\rm L}_\varepsilon$ ) is the
right shift by $\eta$ (resp. the left shift by $\varepsilon$)
in the monoid $M$ (resp. $N$).

The equalities (\ref{6})-(\ref{6a}) together with the fact that
$\varphi$ and $\varphi^{-1}$ are iso's imply

\begin{Prop}
\label{P1}
Both $f$ and $g$ are injective morphisms of monoids, whereas
${\rm R}_\eta$ and ${\rm L}_\varepsilon$ are surjective maps
of sets.
\end{Prop}

Setting now $n = n' = 1$  in (\ref{3})  one gets
\begin{userlabel}
\begin{equation}
\label{7}
(\varphi\comp f)(m) = \eta m
\end{equation}
and, dually,
\begin{equation}
\label{7a}
(\varphi^{-1}\comp g)(n) = n\varepsilon
\end{equation}
\end{userlabel}
which, together with Prop.\ref{P1}, implies that both
${\rm L}_\eta$ and ${\rm R}_\varepsilon$ are {\it injective},
so that if there exists an adjunction with $f$ (resp. $g$) non-iso,
then $M$ (resp. $N$) is, at least, non-commutative monoid; it may not
be a group as well.

Let us now reinterpret $\eta$ and $\varepsilon$ as
natural transformations:

\begin{equation}
\label{8}
\eta\colon\sar 1_M,{},g\comp f,;\; \; \ \varepsilon\colon\sar f\comp g,{},1_N;
\end{equation}
(unit and counit of the adjunction $(f,g,\varphi)$).

Indeed, given two morphisms of monoids $\darr {f_1} {f_2} M {M'}$,
a natural transformation $\mu\colon f_1\arr f_2$ is uniquely determined
by an element $\mu'\in M'$ such that the identity

\begin{equation}
\label{9}
f_2(m)\mu' = \mu'f_1(m)
\end{equation}
holds for every $m\in M$; in this case $\mu$ itself can be identified
with the triple $(f_1,\mu',f_2)$ or, by abuse of notations, with $\mu'$
itself.

But identities (\ref{5})-(\ref{5a}) and (\ref{7})-(\ref{7a}) together give:

\begin{userlabel}
\begin{equation}
\label{10}
((g\comp f)(m))\eta = \eta m\ \ \text{ for any }m\in M
\end{equation}
and, dually,

\begin{equation}
\label{10a}
\varepsilon((f\comp g)(n)) = n\varepsilon\ \ \text{ for any }n\in N
\end{equation}
\end{userlabel}
which exactly states that $\eta$ and $\varepsilon$ define
natural transformations (\ref{8}).

Moreover, one must have the identities:

\begin{userlabel}
\begin{equation}
\label{11}
1_M = g(\varepsilon)\eta
\end{equation}

\begin{equation}
\label{11a}
1_N = \varepsilon f(\eta).
\end{equation}
\end{userlabel}
Finally, the adjunction $(f,g,\varphi)\colon M\rightharpoonup N$ can be
described by the data
$$
(f\colon M\rightarrow N,g\colon N\rightarrow M,\eta\in M,\varepsilon\in N)
$$
satisfying the identities (\ref{10})-(\ref{11a}) (see \cite{ML}, p.81);
$\varphi$ and $\varphi^{-1}$ are then {\it defined} by eqs. (\ref{5}-\ref{5a}).

Define now the category \admon\ such that its objects are just all
adjunctions\linebreak
$(f,g,\eta,\varepsilon)\colon M\rightharpoonup N$ and, given another
adjunction
$(f',g',\eta',\varepsilon')\colon M'\rightharpoonup N'$,
a pair of monoid morphisms
$(l\colon M\rightarrow M',\ r\colon N\rightarrow N')$
is a morphism
$\spar(f,g,\eta,\varepsilon).(l,r).(f',g',\eta',\varepsilon');$
in \admon\ if the diagram

\begin{equation}
\label{12}
\divide\dgARROWLENGTH by2
\begin{diagram}
\node{M}\arrow{e,t}{f}\arrow{s,l}{l}\node{N}\arrow{e,t}{g}\arrow{s,l}{r}%
\node{M}\arrow{s,l}{l}\\
\node{M'}\arrow{e,t}{f'}\node{N'}\arrow{e,t}{g'}\node{M'}
\end{diagram}
\end{equation}
is commutative and, besides, if

\begin{equation}
\label{13}
r(\varepsilon) = \varepsilon',\ \ l(\eta) = \eta'.
\end{equation}

Denote by \mon\ the category of monoids. One immediately sees that there are
two forgetful functors

\begin{equation}
\label{14}
{\rm L},{\rm R}\colon\admon\arr\mon
\end{equation}
defined as follows:

\begin{equation}
\label{15}
{\rm L}(l,r) = l;\ \ \ {\rm R}(l,r) = r.
\end{equation}

Given now an adjunction
$(f,g,\eta,\varepsilon)\colon M\rightharpoonup N$
and considering the commutative diagram
\begin{equation}
\label{16}
\divide\dgARROWLENGTH by2
\begin{diagram}
\node{M}\arrow{e,t}{f}\arrow{s,l}{1_M}\node{N}\arrow{e,t}{g}\arrow{s,l}{I}%
\node{M}\arrow{s,l}{1_M}\\
\node{M}\arrow{e,t}{}\node{{\rm Im}(g)}\arrow{e,t}{\subset}\node{M}
\end{diagram}
\end{equation}
where $I$ is an iso due to Prop.\ref{P1}, one can see that
$(I\comp f,g'=({\rm Im}(g)\subset M),\eta,\varepsilon'=g(\varepsilon))$
is an adjunction isomorphic to $(f,g,\eta,\varepsilon)$ as an object of
\admon, where the isomorphism is given by $(1_M,I)$; that all this iso's
together generate the natural equivalence of \admon\ with its full
subcategory consisting of just those adjunctions
$(f,g,\eta,\varepsilon)\colon M\rightharpoonup N$ in which $N$ is
a submonoid of $M$ and $g$ is the inclusion map $N\subset M$.

Define the category \admonl\ as follows: objects of \admonl\ are all data
of the type
$(N\subset M,f\colon M\rightarrow M, \eta\in M,\varepsilon\in M)$,
where $M$ is a monoid, $N$ its submonoid, $f$ a monoid
endomorphism such that $f(M)\subset N$ and, besides, the identities
\begin{alphalabel}
\label{17}
\begin{eqnarray}
\label{17a}
f(m)\eta & = & \eta m\ \ \ (m\in M)\\
\label{17b}
\varepsilon f(n) & = & n\varepsilon\ \ \ (n\in N)\\
\label{17c}
\varepsilon\eta & = & 1\\
\label{17d}
\varepsilon f(\eta) & = & 1
\end{eqnarray}
\end{alphalabel}
hold. In other words, objects of \admonl\ are just monoids equipped with
some additional structure (a submonoid $N\subset M$, an endomorphism
$f\colon M\arr M$ such that $f(M)\subset N$ and elements $\eta\in M$,
$\varepsilon\in M$ satisfying eqs.(\ref{17a})-(\ref{17d})); then a morphism
in \admonl\ is just a monoid morphism respecting this structure.
So we have (see above):
\begin{Prop}
\label{P2}
The category \admon\ is naturally equivalent to the category \admonl.
\end{Prop}

Note that though the transition from \admon\ to \admonl\ breaks the
``$\Z/2$-symmetry'' ($f\leftrightarrow g$, $\varepsilon\leftrightarrow
\eta$) we get simpler objects instead and simpler relations
(\ref{17a})-(\ref{17d}) instead of ``symmetric'' ones
(\ref{10})-(\ref{11a}).

Now we will give some conditions on an object
$(M\supset N,f,\eta,\varepsilon)$
equivalent to the statement that $f$ is an isomorphism.
\begin{Prop}
\label{P3}
Let
$(M\supset N,f,\eta,\varepsilon)$
be an object of \admonl. Then the following conditions {\rm(a)-(g)}
are equivalent:

\noindent{\rm(a)} $f$ is surjective;

\noindent{\rm(b)} $f$ is an isomorphism;

\noindent{\rm(c)} $N = M$;

\noindent{\rm(d)} $f(\eta) = \eta$;

\noindent{\rm(e)} $f(\varepsilon) = \varepsilon$;

\noindent{\rm(f)} $\eta\varepsilon = 1$
   {\rm(i.e., $\eta$ is invertible in $M$ due to eq.(\ref{17c}))};

\noindent{\rm(g)} for every $m\in M$ one has $f(m) = \eta m\varepsilon$
   {\rm(i.e., $f$ is an inner automorphism of $M$ due to {\rm(f)})}.
\end{Prop}

\noindent{\bf Proof.} (a)$\Longleftrightarrow$(b) due to Prop.\ref{P1};
(b)$\Longrightarrow$(c) is evident, because $N$ contains $f(M) = M$;

(c)$\Longrightarrow$(d): Multiplying eq.(\ref{17a}) by $\varepsilon$
from the left (resp. multiplying eq.(\ref{17b}) by $f(\eta)$) from the right)
one gets:
\begin{alphalabel}
\label{18}
\begin{eqnarray}
\label{18a}
m = \varepsilon f(m)\eta\ \ \ (m\in M),\\
\label{18b}
n = \varepsilon f(n\eta)\ \ \ (n\in N).
\end{eqnarray}
\end{alphalabel}

If $N = M$, then due to eq.(\ref{18b}) every $m\in M$ has a representation
$m = \varepsilon f(m')$ for some $m'\in M$. In particular,
\begin{equation}
\label{19}
\eta = \varepsilon f(\eta^2) = (\varepsilon f(\eta))f(\eta) = f(\eta);
\end{equation}

(d)$\Longrightarrow$(e): $\varepsilon\in N$ implies
\begin{equation}
\label{20}
\varepsilon = \varepsilon1 = \varepsilon f^2(1) \stackrel{(\ref{17c})}{=}
\varepsilon f^2(\varepsilon\eta) = \varepsilon f^2(\varepsilon)f^2(\eta)\\
\stackrel{(d)}{=} \varepsilon f^2(\varepsilon)\eta
\stackrel{(\ref{17a})}{=} \varepsilon\eta f(\varepsilon)
\stackrel{(\ref{17c})}{=} f(\varepsilon);
\end{equation}

(e)$\Longrightarrow$(f): Indeed:
\begin{equation}
\label{21}
\eta\varepsilon\stackrel{(\ref{17a})}{=} f(\varepsilon)\eta
\stackrel{(e)}{=} \varepsilon\eta = 1;
\end{equation}

(f)$\Longrightarrow$(g): One has: $f(m)\stackrel{(f)}{=}
f(m)\eta\varepsilon\stackrel{(\ref{17a})}{=}\eta m\varepsilon$;

(g)$\Longrightarrow$(b), because
$1 = f(1) = \eta\varepsilon$, i.e., (g)$\Longrightarrow$(f) and defining
\begin{equation}
\label{22}
f^{-1}(m) := \varepsilon m\eta
\end{equation}
one sees that
$f^{-1}f(m) = ff^{-1}(m) = m$\ \eop

Return now to the definition of \admonl.  One sees that \admonl\ is
``almost'' the variety of algebras in the Birkhoff's sense (see, e.g.,
\cite{C}).  In more detail, let
$$
\Omega =
\Omega_0\cup\Omega_1\cup\Omega_2, \;\text{where $\Omega_0 =
\{1,\eta,\varepsilon\}$, $\Omega_1 = \{f\}$ and $\Omega_2 = \{\mult\}$
}
$$
i.e, $\Omega_0$, $\Omega_1$ and $\Omega_2$ are, respectively, the set of 0-ary,
unary and binary operations. The
set of ``equations'' consists of equations stating that $\mult$ and
$1$ determine a monoid structure, $f$ is an endomorphism of the
corresponding structure, with eqs.(\ref{17a})-(\ref{17d}) added. One
sees that eq.(\ref{17b}) is not an equation in Birkhoff's sense,
because $n$ there is {\it restricted} to the subset $N\subset M$.
This means that, generally speaking, \admonl\ is a variety of
``sorted'' algebras with two sorts of algebras: ``$N$-like'' and
``$M$-like''.  But, noting that due to eq.(\ref{5a}):
\begin{equation}
\label{23}
N = \varepsilon f(M)
\end{equation}
one can derive from eq.(\ref{17b}) that the equality
\begin{equation}
\label{24}
\varepsilon f(\varepsilon f(m)) = \varepsilon f(m)\varepsilon\ \ \ (m\in M)
\end{equation}
is valid in which $m$ runs over the {\it whole} $M$.

We will prove now that vice versa, the data
$(M,f\colon M\arr M,\eta,\varepsilon)$ together with equations
(\ref{17a}), (\ref{17c})-(\ref{17d}), equations stating that $f$
is an automorphism of monoids as well as eq.(\ref{24})
(instead of eq.(\ref{17b})), reconstruct the remaining data,
namely, the submonoid $N$ and the equation (\ref{17b}) valid on $N$.

Indeed, {\it define} $N$ by eq.(\ref{23}) as a set; we have to prove that
this set is, in fact, a submonoid (eq.(\ref{17b}) for $n\in N$ follows
immediately from eq.(\ref{24})).

We see, first of all, that $1\in N$ due to eq.(\ref{17d}). Suppose now
that $n_1,n_2\in N$, i.e., for some $m_1,m_2\in M$ one has
$n_1 = \varepsilon f(m_1)$, $n_2 = \varepsilon f(m_2)$. Then:
\begin{equation}
\label{25}
n_1n_2 = \varepsilon f(m_1)\varepsilon f(m_2)
\stackrel{(\ref{24})}{=}
\varepsilon f(\varepsilon f(m_1))f(m_2) =
\varepsilon f(\varepsilon f(m_1)m_2)\in N.
\end{equation}

This proves that $N$ is, actually, a submonoid of $M$.

\begin{Prop}
\label{P4}
The category \admon\ is naturally equivalent to the Birkhoff variety
of monoids $M$ equipped with the structure
$(M,\; f\colon M\arr M,\; \eta, \varepsilon\in M)$, where $f$ is an endomorphism
of monoids, satisfying the following conditions:
\begin{alphalabel}
\begin{eqnarray}
\label{26a}
\varepsilon\eta & = & 1,\\
\label{26b}
\varepsilon f(\eta) & = & 1,\\
\label{26c}
\varepsilon f(\varepsilon) & = & \varepsilon^2,\\
\label{26d}
\varepsilon f^2(m) & = & f(m)\varepsilon\ \ \ (m\in M),\\
\label{26e}
f(m)\eta & = & \eta m\ \ \ (m\in M).
\end{eqnarray}
\end{alphalabel}
\end{Prop}
{\bf Proof.} It remains to prove only that eqs.(\ref{26c})-(\ref{26d})
together are equivalent to eq.(\ref{24}) above. Indeed, eq.(\ref{26c})
is a particular case of eq.(\ref{24}) for $m = 1$, whereas eq.(\ref{26d})
is obtained from eq.(\ref{24}) if one substitutes $m = \eta m'$ and takes
into account eq.(\ref{26a}). On the other hand, for $m\in M$ one has:
\begin{equation}
\label{27}
\varepsilon f(\varepsilon f(m)) =
\varepsilon f(\varepsilon)f^2(m)\stackrel{(\ref{26c})}{=}
\varepsilon^2f^2(m)\stackrel{(\ref{26d})}{=}
\varepsilon f(m)\varepsilon\ \eop
\end{equation}

Denote \admonb\ the Birkhoff variety described by Prop.\ref{P4}.
We will prove next that the equality $\eta\varepsilon=1$ is not
satisfied in \admonb. To this end, we will consider in details
the ``minimal model of the theory \admonb'', in other words,
the free \admonb\ algebra $\script F(\emptyset)$ (which is an 
initial object in \admonb).

Define elements $\eta_k,\varepsilon_k\in\script F(\emptyset)$ ($k\in \N$)
as follows:
\begin{alphalabel}
\label{28}
\begin{eqnarray}
\label{28a}
\eta_0 := \eta,\ \ \ \eta_{k+1} := f(\eta_{k})\\
\label{28b}
\varepsilon_0 := \varepsilon,\ \ \ \varepsilon_{k+1} := f(\varepsilon_{k})
\end{eqnarray}
\end{alphalabel}

It is clear that, as monoid, \free\  is generated by elements
$\eta_k,\ \varepsilon_k$.
In other words, the monoid \free\ can be represented as
\begin{equation}
\label{29}
\free = \freemon/{\rm R}
\end{equation}
for some set of relations ${\rm R}$, where
$
\script F_{\bf Mon}(\eta_0,\eta_1,\dots,\eta_k,\dots
\varepsilon_0,\varepsilon_1,\dots,\varepsilon_k,\dots)
$
is the free monoid generated by $2\times\N$ variables $\fam{\eta_k,\varepsilon_m}{k,m\in\N}$.
The following  proposition describes the corresponding set of relations $\rm R$.
\begin{Prop}
\label{P5}
The set of relations ${\rm R}$ in \free\ is generated by the following set
${\rm R}_0$ of relations:
\begin{alphalabel}
\label{30}
\begin{eqnarray}
\label{30a}
\varepsilon_i\varepsilon_j & = & \varepsilon_{j-1}\varepsilon_i\ (j > i),\\
\label{30b}
\eta_j\eta_i & = & \eta_i\eta_{j-1}\ (j > i),\\
\label{30c}
\varepsilon_i\eta_j & = & \left\{
\begin{array}{ll}
\eta_{j-1}\varepsilon_i & (j > i+1)\\
\eta_{j}\varepsilon_{i-1} & (i > j)\\
\varepsilon_{i-1}\eta_i & (i=j > 0)\\
\varepsilon_i\eta_i & (j=i+1)\\
1 & (i = j = 0)
\end{array}
\right.
\end{eqnarray}
\end{alphalabel}
\end{Prop}

{\bf Proof.} One easily checks that eqs.(\ref{30a})-(\ref{30c})
are satisfied, being either particular cases of some of the relations
(\ref{26a})-(\ref{26e}), or can be obtained from the latter ones after
applying $f^k$ to both sides of (\ref{26a})-(\ref{26e}) (for some $k>0$).

Vice versa, define the elements $\eta$ and $\varepsilon$ and
an automorphism $f$ of the monoid \freemon\ by ``inverting''
definitions~(\ref{28}) above:
\begin{alphalabel}
\label{31}
\begin{eqnarray}
\label{31a}
\eta := \eta_0,\ \ \ f(\eta_{k}) := \eta_{k+1,}\\
\label{31b}
\varepsilon := \varepsilon_0,\ \ \ f(\varepsilon_{k}) := \varepsilon_{k+1}.
\end{eqnarray}
\end{alphalabel}
One easily checks, that the automorphism $f$ ``survives'' the factorization
by relations (\ref{30}) above and induction on the length of words proves
that, in  the monoid
\begin{equation}
\label{32}
\freemon/{\rm R}_0,
\end{equation}
relations (\ref{26d})-(\ref{26e}) are satisfied. \eop

From now on we will identify \free\ with the monoid (\ref{32}) equipped
with $\eta$, $\varepsilon$ and $f$ defined by eqs.(\ref{31}).

Returning now to our original problem: one sees that it is exactly equivalent
to the question whether or not the identity
\begin{equation}
\label{33}
\eta\varepsilon = 1
\end{equation}
holds in the monoid \free.

The following theorem provides us with the canonical form for elements
of \free\ and simultaneously gives the negative answer to the last question.

\begin{Th}
\label{T1}
For every element $m$ of \free, there exist the only pair $k,l\in\N$
and the only pair of sequences
\begin{alphalabel}
\begin{equation}
\label{34a}
0\le i_1\le\dots\le i_k,\ \ \ j_1\ge\dots\ge j_l\ge0
\end{equation}
such that
\begin{equation}
\label{34b}
m = \eta_{i_1}\dots\eta_{i_k}\varepsilon_{j_1}\dots\varepsilon_{j_l}.
\end{equation}
\end{alphalabel}
(We asume that if $k=0$ (resp. $l=0$), then the corresponding sequence
in {\rm (\ref{34a})} above is empty and 
the corresponding product in
{\rm (\ref{34b})} is replaced with the neutral element $1$ of the monoid \free.)
\end{Th}

\noindent{\bf Proof.} On \freemon, define the binary relation $\models$
as follows. First of all, set:
\begin{alphalabel}
\label{35}
\begin{eqnarray}
\label{35a}
\varepsilon_i\varepsilon_j & \models & \varepsilon_{j-1}\varepsilon_i\ (j > i)\\
\label{35b}
\eta_j\eta_i & \models & \eta_i\eta_{j-1}\ (j > i)\\
\label{35c}
\varepsilon_i\eta_j & \models & \left\{
\begin{array}{ll}
\eta_{j-1}\varepsilon_i & (j > i+1)\\
\eta_{j}\varepsilon_{i-1} & (i > j)\\
\varepsilon_{i-1}\eta_i & (i=j > 0)\\
\varepsilon_i\eta_i & (j=i+1)\\
1 & (i = j = 0)
\end{array}
\right.
\end{eqnarray}
\end{alphalabel}

Let now $\ge$ be the smallest preorder on \freemon\ containing $\models$
and turning \freemon\ into a preordered monoid.

One sees that $m\ge m'$ if and only if there exists a sequence
\begin{equation}
\label{36}
m = m_0,m_1,\dots,m_n=m'
\end{equation}
such that, for any $0\le i<n$, there exist
$$L,R,\mu,\mu'\in\freemon$$ for which $\mu\models\mu'$ and both
$m_i=L\mu R$ and $m_{i+1}=L\mu'R$.

Let now $\sim$ be the equivalence relation on \freemon\ generated by
the relation $\ge$ (i.e., $m\sim m'$ if and only if both $m$ and $m'$ belong to
the same connected component of the preorder relation $\ge$).

It is rather clear that the relation $\sim$ coincides with the
relation $\rm R$ from eq.(\ref{29}) defining \free.

It is also clear that the r.h.s. of the canonical representation
(\ref{34b}), considered as an element of \freemon\ is a minimal element
with respect to the preorder relation $\ge$. The only thing to prove
is that {\it every equivalence class of the relation $\sim$ contains
the only minimal element with respect to $\ge$.}

To prove this, it suffices to prove that the relation $\ge$
satisfies conditions of reduction theorem of M.~Newman \cite{N},
or its weaker version given in \cite{C}. The latter conditions
on $\ge$ are the following conditions (A) and (B): 
\medskip

\noindent (A) {\it For any $m\in\freemon$, there exists $k\in\N$ such that
for any decreasing sequence
$m=m_0>m_1>\dots>m_{k'}$ one has $k\le k'$};
\medskip

\noindent (B) {\it Any pair $m,m'\in\freemon$ with a common parent is bounded
from below, i.e., there exists $b\in\freemon$ such that both $m\ge b$
and $m'\ge b$\/}.
Here $p$ is said to be a {\it parent} of $m$ if it is the smallest element
such that $p\ge m,\ p\not=m$.
\medskip

To prove (A), consider the morphism of monoids
\begin{alphalabel}
\label{37}
\begin{equation}
\label{37a}
d\colon\freemon\arr\N
\end{equation}
uniquelly determined by the images
\begin{equation}
\label{37b}
d(\eta_i) = i+1,\ \ \ d(\varepsilon_i) = i+1.
\end{equation}
\end{alphalabel}

One easily sees from relations (\ref{35}) that $d$ is a morphism of
preordered monoids (i.e., respects preorders).  Moreover, it is clear
now that the relation $\ge$ is, in fact, an order relation (because $d(m)>d(m')$
for any pair $m,m'\in\freemon$ such that $m\models m'$).  (A) is obvious now.

To prove (B), observe first that
\medskip

\noindent(C) {\it $p$ is a parent of $m$ if and only if $p=L\mu R$, $m=L\mu'R$ and
$\mu\models\mu'$}.
\medskip

Indeed, the l.h.s. of any of the particular cases (\ref{35}) of the relation
$\models$, except, perhaps, the last one, is a parent of its r.h.s. just
because $d({\rm l.h.s.})-d({\rm r.h.s.})=1$; as to the last particular case,
$\varepsilon_0\eta_0=1$, this is the only case such that the length of l.h.s.
$\not=$ the length of r.h.s, which implies that in this case as well the l.h.s.
is the parent of the r.h.s.

Let $m, m'\in\freemon$ have a common parent, $p$.

Now follows the most boring part of all this mess.

There are five possible cases:
\medskip

\noindent(I) $p = a\mu b\nu c$, $m=a\mu'b\nu c$, $m'=a\mu b\nu' c$, where
$\mu\models\mu'$ and $\nu\models\nu'$.
\medskip

Clearly, in this case 
$b:=a\mu'b\nu'c$ is a common lower bound for
both $m$ and $m'$.

In all of the cases (II)-(V) below, $p$ is of the form $L\mu_1\mu_2\mu_3R$,
where every $\mu_i$ is of length 1 (i.e., is either $\eta_j$ for some $j$, or
$\varepsilon_j'$ for some $j'$). In what follows, terms $L$ and $R$ will be
omitted, because they take no explicit part in the process of finding of the
lower bound $b$.

\medskip
\noindent(II) $p=\varepsilon_i\varepsilon_j\varepsilon_k\ \ (i<j<k)$;
\medskip
Applying relation (\ref{35a}) one gets:
\begin{eqnarray*}
p\models m &:=& \varepsilon_{j-1}\varepsilon_i\varepsilon_k \ge
\varepsilon_{j-1}\varepsilon_{k-1}\varepsilon_i \ge
\varepsilon_{k-2}\varepsilon_{j-1}\varepsilon_i\\
p\models m' &:=& \varepsilon_i\varepsilon_{k-1}\varepsilon_j \ge
\varepsilon_{k-2}\varepsilon_i\varepsilon_j \ge
\varepsilon_{k-2}\varepsilon_{j-1}\varepsilon_i
\end{eqnarray*}
So in this case
$b=\varepsilon_{k-2}\varepsilon_{j-1}\varepsilon_i$
is the lower bound of both $m$ and $m'$.

\medskip
\noindent(III) $p=\eta_i\eta_j\eta_k\ \ (i>j>k)$;
\medskip

\noindent This case is  ``dual'' in an obvious sense to case (II).

\medskip
\noindent(IV) $p=\varepsilon_i\varepsilon_j\eta_k\ \ (i<j)$;
\medskip

\noindent This case is subdivided into following 7 subcases below:

\noindent a) $k>j+1$
\begin{eqnarray*}
p\models m &:=& \varepsilon_{j-1}\varepsilon_i\eta_k\ge
\varepsilon_{j-1}\eta_{k-1}\varepsilon_i\ge
\eta_{k-2}\varepsilon_{j-1}\varepsilon_i\\
p\models m' &:=& \varepsilon_i\eta_{k-1}\varepsilon_j\ge
\eta_{k-2}\varepsilon_i\varepsilon_j\ge
\eta_{k-2}\varepsilon_{j-1}\varepsilon_i
\end{eqnarray*}
i.e.,
$
b=\eta_{k-2}\varepsilon_{j-1}\varepsilon_i
$
is a common lower bound of both $m$ and $m'$.

\noindent b) $k=j+1$
\begin{eqnarray*}
p=\varepsilon_i\varepsilon_j\varepsilon_{j+1}\models m &:=&
\varepsilon_{j-1}\varepsilon_i\eta_{j+1}\ge
\varepsilon_{j-1}\eta_j\varepsilon_i\ge\dots\ge\varepsilon_i\\
p\models m' &:=& \varepsilon_i\varepsilon_j\eta_j\ge\dots\ge\varepsilon_i
\end{eqnarray*}
i.e.,
$
b=\varepsilon_i
$
is a common lower bound of both $m$ and $m'$.

\noindent c) $k=j$
\begin{eqnarray*}
p=\varepsilon_i\varepsilon_j\eta_j\models m &:=&
\varepsilon_{j-1}\varepsilon_i\eta_j\ge
\left\{
\begin{array}{ll}
\varepsilon_{j-1}=\varepsilon_i & (i = j-1)\\
\varepsilon_{j-1}\eta_{j-1}\varepsilon_i\ge\dots\ge\varepsilon_i & (i<j-1)
\end{array}
\right.\\
p\models m' &:=& \varepsilon_i\varepsilon_{j-1}\eta_j\ge\dots\ge\varepsilon_i
\end{eqnarray*}
i.e.,
$
b=\varepsilon_i
$
is a common lower bound of both $m$ and $m'$.

\noindent d) $i+1<k<j$
\begin{eqnarray*}
p\models m &:=& \varepsilon_{j-1}\varepsilon_i\eta_k\ge
\varepsilon_{j-1}\eta_{k-1}\varepsilon_i\ge
\eta_{k-1}\varepsilon_{j-2}\varepsilon_i\\
p\models m' &:=& \varepsilon_i\eta_k\varepsilon_{j-1}\ge
\eta_{k-1}\varepsilon_{j-2}\varepsilon_i
\end{eqnarray*}
i.e.,
$
b=\eta_{k-1}\varepsilon_{j-2}\varepsilon_i
$
is a common lower bound of both $m$ and $m'$.

\noindent e) $k=i+1<j$
\begin{eqnarray*}
p=\varepsilon_i\varepsilon_j\eta_{i+1}\models m &:=&
\varepsilon_{j-1}\varepsilon_i\eta_{i+1}\ge\dots\ge\varepsilon_{j-1}\\
p\models m' &:=& \varepsilon_i\eta_{i+1}\varepsilon_{j-1}\ge\dots\ge
\varepsilon_{j-1}
\end{eqnarray*}
i.e.,
$
b=\varepsilon_{j-1}
$
is a common lower bound of both $m$ and $m'$.

\noindent f) $k=i$
\begin{eqnarray*}
p\models m &:=& \varepsilon_{j-1}\varepsilon_i\eta_i\ge\dots\ge
\varepsilon_{j-1}\\
p\models m' &:=& \varepsilon_i\eta_i\varepsilon_{j-1}\ge\dots\ge
\varepsilon_{j-1}\\
\end{eqnarray*}
i.e.,
$
b=\varepsilon_{j-1}
$
is a common lower bound of both $m$ and $m'$.

\noindent g) $k<i$
\begin{eqnarray*}
p\models m &:=& \varepsilon_{j-1}\varepsilon_i\eta_k\ge
\varepsilon_{j-1}\eta_k\varepsilon_{i-1}\ge
\eta_k\varepsilon_{j-2}\varepsilon_{i-1}\\
p\models m' &:=& \varepsilon_i\eta_k\varepsilon_{j-1}\ge
\eta_k\varepsilon_{j-2}\varepsilon_{i-1}
\end{eqnarray*}
i.e.,
$
\eta_k\varepsilon_{j-2}\varepsilon_{i-1}
$
is a common lower bound of both $m$ and $m'$.

\medskip
\noindent(V) $p=\varepsilon_i\eta_j\eta_k\ \ (j>k)$;
\medskip

\noindent This case is in a sense dual to case (IV) and is
subdivided into following 7 subcases below:

\noindent a) $i>j$
\begin{eqnarray*}
p\models m &:=& \varepsilon_i\eta_k\eta_{j-1}\ge
\eta_k\varepsilon_{i-1}\eta_{j-1}\ge
\eta_k\eta_{j-1}\varepsilon_{i-2}\\
p\models m' &:=& \eta_j\varepsilon_{i-1}\eta_k\ge
\eta_j\eta_k\varepsilon_{i-2}\ge
\eta_k\eta_{j-1}\varepsilon_{i-2}
\end{eqnarray*}
i.e.,
$
b=\eta_k\eta_{j-1}\varepsilon_{i-2}
$
is a common lower bound of both $m$ and $m'$.

\noindent b) $i=j$
\begin{eqnarray*}
p=\eta_j\eta_j\varepsilon_k\models m &:=&
\varepsilon_{j-1}\eta_j\eta_k\ge\dots\ge\eta_k\\
p\models m' &:=& \varepsilon_j\eta_k\eta_{j-1}\ge
\eta_k\varepsilon_{j-1}\eta_{j-1}\ge\dots\ge\eta_k
\end{eqnarray*}
i.e.,
$
b=\eta_k
$
is a common lower bound of both $m$ and $m'$.

\noindent c) $i=j-1$
\begin{eqnarray*}
p=\varepsilon_{j-1}\eta_j\eta_k\models m &:=&
\varepsilon_{j-1}\eta_{j-1}\eta_k\ge\dots\ge\eta_k\\
p\models m' &:=&
\varepsilon_{j-1}\eta_k\eta_{j-1}\ge
\left\{
\begin{array}{ll}
\eta_{j-1}=\eta_k & (k=j-1)\\
\eta_k\varepsilon_{j-2}\eta_{j-1}\ge\dots\ge\eta_k & (k<j-1)
\end{array}
\right.
\end{eqnarray*}
i.e.,
$
b=\eta_k
$
is a common lower bound of both $m$ and $m'$.

\noindent d) $k<i<j-1$
\begin{eqnarray*}
p\models m &:=& \varepsilon_i\eta_k\eta_{j-1}\ge
\eta_k\varepsilon_{i-1}\eta_{j-1}\ge
\eta_k\eta_{j-2}\varepsilon_{i-1}\\
p\models m' &:=& \eta_{j-1}\varepsilon_i\eta_k\ge
\eta_{j-1}\eta_k\varepsilon_{i-1}\ge
\eta_k\eta_{j-2}\varepsilon_{i-1}
\end{eqnarray*}
i.e.,
$
b=\eta_k\eta_{j-2}\varepsilon_{i-1}
$
is a common lower bound of both $m$ and $m'$.

\noindent e) $k=i<j-1$
\begin{eqnarray*}
p=\varepsilon_i\eta_j\eta_i\models m &:=&
\varepsilon_i\eta_i\eta_{j-1}\ge\dots\ge\eta_{j-1}\\
p\models m' &:=& \eta_{j-1}\varepsilon_i\eta_i\ge\dots\ge\eta_{j-1}
\end{eqnarray*}
i.e.,
$
b=\eta_{j-1}
$
is a common lower bound of both $m$ and $m'$.

\noindent f) $k=i+1$
\begin{eqnarray*}
p=\varepsilon_i\eta_j\eta_{i+1}\models m &:=&
\varepsilon_i\eta_{i+1}\eta_{j-1}\ge\dots\ge
\eta_{j-1}\\
p\models m' &:=& \eta_{j-1}\varepsilon_i\eta_{i+1}\ge\dots\ge
\eta_{j-1}\\
\end{eqnarray*}
i.e.,
$
b=\eta_{j-1}
$
is a common lower bound of both $m$ and $m'$.

\noindent g) $k>i+1$
\begin{eqnarray*}
p\models m &:=& \varepsilon_i\eta_k\eta_{j-1}\ge
\eta_{k-1}\varepsilon_i\eta_{j-1}\ge
\eta_{k-1}\eta_{j-2}\varepsilon_i\\
p\models m' &:=& \eta_{j-1}\varepsilon_i\eta_k\ge
\eta_{j-1}\eta_{k-1}\varepsilon_i\ge
\eta_{k-1}\eta_{j-2}\varepsilon_i
\end{eqnarray*}
i.e.,
$
\varepsilon_k\eta_{j-2}\eta_{i-1}
$
is a common lower bound of both $m$ and $m'$.
\eop

\enddocument

\begin{thebibliography}{aaa}
\bibitem{M} Ernest G. Manes {\sl Algebraic Theories}, Springer, (1976)
\bibitem{ML} S. Mac Lain {\sl Categories for the Working Mathematicians},
Springer (1971);
\bibitem{C} P. M. Cohn {\sl Universal Algebra}, Harper's series in Modern Math.,
Harper \& Row (1965);
\bibitem{N} M. H. A.Newman {\sl On Theories with a Combinatorial Definition
of ``Equivalence''}, Ann. of Math., {\bf 43}, p.223--243 (1942).
\end{thebibliography}
\end{document}